%% file: Bianchi.tex
\documentclass[11pt]{article}

\input{macros_angl.tex}

\title{On the representation of integers by 
indefinite binary Hermitian forms} 
\author{Jouni Parkkonen \and Fr\'ed\'eric Paulin}
\date{\today}

\begin{document}
\bibliographystyle{../alphanum}

\maketitle

\begin{abstract} 
Given an integral indefinite binary Hermitian form $f$ over an
  imaginary quadratic number field, we give a precise asymptotic
  equivalent to the number of
  nonequivalent representations, satisfying some congruence properties,
  of the rational integers with absolute value at most $s$ by
  $f$, as $s$ tends to $+\infty$.
\end{abstract}

\section{Introduction}
\label{sect:intro}

Though less thoroughly developped than the real case of binary
quadratic forms initiated by Gauss, the problem of the representation
of integers by integral binary Hermitian forms, along with their
reduction theory, initiated by Hermite, Bianchi and especially
Humbert, has been much studied (see for instance
\cite[Sect.~9]{ElsGruMen98} and references therein). The average, over
the representatives of the binary Hermitian forms with a given
discriminant, of the number of their nonequivalent representations of a
given integer has been computed by Elstrodt, Grunewald, Mennicke
\cite{ElsGruMen87}. In this paper, we concentrate on a given form, and
our result gives a precise asymptotic on the number of nonequivalent
proper representations of rational integers with absolute value at
most $s$ by a given integral indefinite binary Hermitian form.

A binary Hermitian form naturally gives rise to a quaternary quadratic
form. The representations of integers by positive definite quaternary
quadratic forms have been studied for a long time (including
Lagrange's four square theorem, see also the work of Ramanujan as in
\cite{Kloosterman27}.)  In the case of indefinite forms, the counting
problem is complicated by the presence of an infinite group of
automorphs of the form.  General formulas are known (by Siegel's mass
formula, see for instance \cite{EskRudSar91}), but it does not seem to
be easy (or even doable) to deduce our asymptotic formulae from
them. Our proof is geometric, while the methods of \cite{ElsGruMen87}
for the result quoted above rely on analytic techniques for the
appropriate counting Dirichlet series. There are numerous results on
counting integer points with bounded norm on quadrics, see for
instance \cite{deLury38,DukeRudSar93,EskMcMul93,BorovoiRudnick95} and
the excellent survey \cite{Babillot02a}. In this paper, we consider
a problem of a somewhat different nature, and we count appropriate
orbits of integer points on which the form is constant.

Let $K$ be an imaginary quadratic number field, with discriminant
$D_K$, ring of integers $\OOO_K$, and Dedekind zeta function
$\zeta_K$. Fix an indefinite binary Hermitian form $f:\CC^2\ra\RR$
with
\begin{equation}\label{eq:form}
f(u,v)=a|u|^2+2\,\Re(bu\ov v) +c|v|^2 
\end{equation}
which is integral over $K$ (its coefficients satisfy $a,b\in\ZZ$ and
$b\in\OOO_K$).  The discriminant $\Delta(f)=|b|^2-ac$ of the form $f$ is
positive.  The group $\autom$ of automorphs of $f$ consists of those
elements $g\in\SLOK$ for which $f\circ g=f$.

Let $\P_K$ be the set of relatively prime pairs of integers of
$K$. For every $s>0$, the number of nonequivalent proper
representations of rational integers with absolute value at most $s$
by $f$ is
$$
\psi_f(s)=\card\;\;_{\mbox{$\autom$}}\bs
\big\{(u,v)\in\P_K\;:\;|f(u,v)|\leq s\big\}\;.
$$
The finiteness of $\psi_f(s)$ follows from general results on orbits
of algebraic groups defined over number fields
\cite[Lem.~5.3]{BorelHarishChandra62}, see also
\cite[Theo.~10.3]{Shimura04}.

In this paper, we prove the following theorem and its generalization
(Corollary \ref{coro:maingeneral}) for subsets of $\P_K$ satisfying
additional congruence properties.

\btheo \label{theo:mainHerm}
As $s$ tends to $+\infty$, we have the equivalent
$$
\psi_f(s)\sim \frac{\pi\;\covol(\autom)}
{2\;|D_K|\;\zeta_K(2)\;\Delta(f)}\;\;\;s^2\;.
$$
\etheo

Note that the image of $\autom$ in $\PSLC$ is an arithmetic Fuchsian
subgroup, and by definition, $\covol(\autom)$ is the area of the
quotient of the real hyperbolic plane $\C$ with constant curvature
$-1$ preserved by $\autom$.  The main input to prove Theorem
\ref{theo:mainHerm} is the work \cite{ParPau} (building on
\cite{EskMcMul93}), where we proved an equidistribution result for the
boundary of big tubular neighbourhoods of a finite volume totally
geodesic submanifold (here the image of $\C$) in the quotient of a
real hyperbolic space by a lattice (here the Bianchi group $\PSLOK$).

The covolume of the group of automorphs could be computed using
Prasad's very general formula in \cite{Prasad89}. Using the work of
Maclachlan and Reid \cite{Maclachlan84,MacRei91,MacRei03}, building on
results of Humbert and Vigneras, we give an expression for
$\covol(\autom)$ at the end of Section \ref{sect:rephermit} (Remark 1).

As the final result of the note, we indicate how the results of
\cite{MacRei91} can be used to obtain an even more precise expression
of the asymptotic formula in Theorem \ref{theo:mainHerm} when $K=
\QQ(i)$.  A constant $\const(f)\in\{1,2,3,6\}$ is defined as follows.
If $\discr(f)\equiv 0\mod 4$, let $\const(f)=2$. If the coefficients
$a$ and $c$ of the form $f$ as in Equation \eqref{eq:form} are both
even, let $\const(f)=3$ if $\discr(f)\equiv 1\mod 4$, and let
$\const(f)$ be the remainder modulo $8$ of $\discr(f)$ if
$\discr(f)\equiv 2 \mod 4$.
In all other cases, let $\const(f)=1$.

\bcoro \label{coro:gaussian} Let $f$ be an indefinite binary Hermitian
form with Gaussian integral coefficients. Then as $s$ tends to
$+\infty$,
$$
\psi_f(s)\sim\frac{\pi^2}{8\;\const(f)\;\zeta_{\QQ(i)}(2)}\,
\prod_{p|\discr(f)} \big(1+\bigg(\frac{-1}p\bigg)p^{-1}\big)\;\;s^2\;,
$$
where $p$ ranges over the odd positive rational primes and
$\big(\frac{-1}p\big)$ is a Legendre symbol.  
\ecoro

\section{Representing integers by indefinite binary 
Hermitian forms}
\label{sect:rephermit}

Let $K,D_K,\OOO_K,\zeta_K,\P_K$ be as in the introduction, and let
$\omega_K$ be the number of roots of unity in $\OOO_K$, that is,
$\omega_{K}=4$ if $D_K=-4$, $\omega_{K}=6$ if $D_K=-3$ and
$\omega_K=2$ if $ D_K\neq -3,-4$. Let us first recall some facts about
binary Hermitian forms.

The Lie group $\SLC$ acts linearly on the left on $\CC^2$, and it acts
on the right on the set of binary Hermitian forms $f$ by
precomposition, that is by $f\mapsto f\circ g$ for every $g\in \SLC$,
preserving the discriminant: $\Delta(f\circ g)=\Delta(f)$ for every
$g\in \SLC$.  The (nonuniform) lattice $\GaK=\SLOK$ of $\SLC$
preserves the subset $\P_K$ of $\CC^2$ as well as the set of integral
indefinite binary Hermitian forms over $K$. The stabilizer in $\GaK$
of such a form $f$ is the group of automorphs $\autom$ defined in the
introduction.  

For every indefinite binary Hermitian form $f$ as in Equation
\eqref{eq:form} with discriminant $\Delta=\Delta(f)$, 
let
$$
\C_\infty(f)=\{[u:v]\in\PP^1(\CC)\;:\;f(u,v)=0\}
\vspace{-0.4cm}
$$%
\vspace{-0.2cm}
and
$$
\C(f)=\{(z,t)\in\CC\times\;]0,+\infty[\;:\;f(z,1)+|a|\,t^2=0\}\;.
$$
Identifying, as usual, $\PP^1(\CC)$ with $\CC\cup\{\infty\}$, the set
$\C_\infty(f)$ is the circle of center $-\frac{\overline{b}}{a}$ and
radius $\frac{\sqrt{\Delta}}{|a|}$ if $a\neq 0$, and $\C_\infty(f)$ is
the union of a real line with $\{\infty\}$ otherwise. The map
$f\mapsto \C_\infty(f)$ induces a bijection between the set of
indefinite binary Hermitian forms up to multiplication by a nonzero
real factor and the set of circles and real lines in
$\CC\cup\{\infty\}$.  The linear action of $\SLC$ on $\CC^2$ induces a
left action of $\SLC$ by homographies on the set of circles and real
lines in $\PP^1(\CC)$, and the map $f\mapsto \C_\infty(f)$ is
anti-equivariant for the two actions of $\SLC$, in the sense that, for
every $g\in\SLC$,
\begin{equation}\label{eq:antiequiv}
\C_\infty(f\circ g)= g^{-1}\,\C_\infty(f)\;.
\end{equation}

Given a finite index subgroup $G$ of $\GaK$, an integral binary
Hermitian form $f$ is called {\it $G$-reciprocal} if there exists an
element $g$ in $G$ such that $f\circ g=-f$. We define $R_G(f)=2$ if
$f$ is $G$-reciprocal, and $R_G(f)=1$ otherwise.  The values of
$f(z,1)$ are positive on one of the two components of
$\PP^1(\CC)-\C_\infty(f)$ and negative on the other. As the
signs are switched by precomposition by an element $g$ as above, the
reciprocity of the form $f$ is equivalent to saying that there exists
an element of $G$ preserving $\C_\infty(f)$ and exchanging the two
complementary components of $\C_\infty(f)$.

\bigskip Theorem \ref{theo:mainHerm} will follow from the following
more general result. For every finite index subgroup $G$ of $\Ga_K$,
let $G_{(1,0)}$ be the stabilizer of $(1,0)\in\CC^2$ in $G$; let
$\iota_G=1$ if $-\id\in G$, and $\iota_G=2$ otherwise; for every
$s>0$, let
$$
\psi_{f,G}(s)=\card\;\;_{\mbox{$\autom\cap G$}}\bs
\big\{(u,v)\in G(0,1)\;:\;|f(u,v)|\leq s\big\}\;.
$$
Taking $G=\GaK$, which acts transitively on $\P_K$, we obtain the
counting function $\psi_f=\psi_{f,\GaK}$ of the introduction. Note that
the image of $\autom\cap G$ in $\PSLC$ is again an arithmetic Fuchsian
subgroup. Let $\Ga_{K,\infty}$ be the stabilizer of $[1:0]\in
\PP^1(\CC)$ in $\GaK$.

\btheo\label{theo:mainversionG} Let $f$ be an integral indefinite
binary Hermitian form over an imaginary quadratic number field $K$, and
let $G$ be a finite index subgroup of $\SLOK$. Then as $s$ tends to
$+\infty$, we have the equivalent
$$
\psi_{f,G}(s)\sim 
\frac{\pi\;\iota_G\;[\Ga_{K,\infty}:G_{(1,0)}]\;\covol(\autom\cap G)}
{2\;\omega_K\;|D_K|\;\zeta_K(2)\;\Delta(f)\;[\GaK:G]}\;\;\;s^2\;.
$$
\etheo

\dem Let us first recall a geometric result from \cite{ParPau} that
will be used to prove this theorem.  A subset $A$ of a set endowed
with an action of a group $G$ is {\it precisely invariant} under this
group if for every $g\in G$, if $gA\cap A$ is nonempty, then $gA=A$.

Let $n\geq 2$ and let $\hnr$ be the upper halfspace model of the real
hyperbolic space of dimension $n$, with (constant) sectional curvature
$- 1$. Let $F$ be a finite covolume discrete group of isometries of
$\hnr$. Let $1\leq k\leq n-1$ and let $\C$ be a real hyperbolic
subspace of dimension $k$ of $\hnr$, whose stabilizer $F_\C$ in $F$
has finite covolume. Let $\H$ be a horoball in $\hnr$, which is
precisely invariant under $F$, with stabilizer $F_\H$.

For every $\alpha,\beta\in F$, denote by $\delta_{\alpha,\beta}$ the
common perpendicular arc between $\alpha \C$ and the horosphere
$\beta\partial\H$, and let $\ell(\delta_{\alpha,\beta})$ be its
length, counted positively if $\delta_{\alpha,\beta}$ exits $\beta\H$
at its endpoint on $\beta\partial\H$, and negatively otherwise. By
convention, $\ell(\delta_{\alpha,\beta})=-\infty$ if the boundary at
infinity of $\alpha \C$ contains the point at infinity of $\beta\H$.
Also define the multiplicity of $\delta_{\alpha,\beta}$ as $m(\alpha,
\beta) = 1/{\card(\alpha\,F_\C\alpha^{-1}\; \cap\; \beta\,
  F_\H\beta^{-1})}$.  For every $t\geq 0$, define $\N(t)=
\N_{F,\C,\H}(t)$ as the number, counted with multiplicity, of the
orbits under $F$ in the set of the common perpendicular arcs
$\delta_{\alpha,\beta}$ for $\alpha,\beta\in F$ with length at most
$t$.

For every $m\in\NN$, denoting by  $\SSS_m$  the unit
sphere of the Euclidean space $\RR^{m+1}$, endowed with its induced
Riemannian metric, we have the following result:

\btheo [{\cite[Coro.~4.9]{ParPau}}]\label{theo:mainBPP}
As $t\to+\infty$, we have
$$
\N(t)\sim\frac{\Vol(\SSS_{n-k-1})\Vol(F_\H\bs\H)\Vol(F_\C\bs\C)}
{\Vol(\SSS_{n-1})\Vol(F\bs\hnr)}\;e^{(n-1)t}\;.\;\;\Box
$$
\etheo

\bigskip Let $f,K$ and $G$ be as in the statement of Theorem
\ref{theo:mainversionG}. We write $f$ as in Equation \eqref{eq:form},
and denote by $\Delta$ its discriminant. In order to apply Theorem
\ref{theo:mainBPP}, we first define the various objects $n$, $k$, $F$,
$\H$, and $\C$ that appear in its statement.

Let $n=3$ and $k=2$, so that $\Vol(\SSS_{n-1})=4\pi$,
$\Vol(\SSS_{n-k-1}) = 2$,  the boundary at infinity of $\hnr$ is
$\phtr=\CC\cup\{\infty\}$, and $\PSLC$  acts faithfully and isometrically
on $\hnr$ by the Poincar\'e extension of homographies. 

For any subgroup $H$ of $\SLC$, we denote by $\pr H$ its image in
$\PSLC$, except that the image of $\autom$ is denoted by
$\Pautom$. Let $F=\pr G$.

The {\it Bianchi group} $\pr{\Ga_K}= \PSLOK$ acts discretely on
$\htr$, with finite covolume. By a formula essentially due to Humbert
(see for instance sections 8.8 and 9.6 of \cite{ElsGruMen98}), we have
$$
\Vol(\pr{\Ga_K}\,\bs\htr)=\frac{1}{4\pi^2}|D_K|^{3/2}\zeta_K(2)\;.
$$
Note that $\Vol(F\bs\htr)=[\pr{\Ga_K}:F]\Vol(\pr{\Ga_K}\,\bs\htr)$ and
$[\pr{\Ga_K}:\pr G]= \frac{1}{\iota_G}[\Ga_K:G]$ by the definition of
$\iota_G$.

Let $\H$ be the horoball centered at $\infty$ in $\htr$ that consists
of all points in $\htr$ of Euclidean height at least $1$. This
horoball is precisely invariant under $\pr{\Ga_K}$ by Shimizu's lemma
\cite{Shimizu63}. The stabilizer in the Bianchi group of the point at
infinity $\infty$ is equal to $(\pr{\Ga_K})_\H$. A fundamental domain
$\Omega_K$ for the action of $(\pr{\Ga_K})_\H$ on $\CC$ is
\begin{itemize}
 \item the rectangle $[0,1]\times
[0,\frac{\sqrt{|D_K|}}{2}]$ if $D_K\neq-3,-4$, 
\item the rectangle
$[0,1]\times[0,\frac{1}{2}]$ if $D_K=-4$, 
\item and the parallelogram with
vertices $0,\frac{1}{2}-\frac{i}{2\sqrt{3}},
\frac{1}{2}+\frac{i}{2\sqrt{3}}, \frac{i}{\sqrt{3}}$ if $D_K= -3$, 
\end{itemize}
see for instance \cite[page 318]{ElsGruMen98}. Since $\F_K= \Omega_K
\times [1,+\infty[\;\subset \htr$ is a fundamental domain for the
action of $(\pr{\Ga_K})_\H$ on $\H$, an easy volume computation in
hyperbolic geometry gives $\Vol(\F_K)= \frac{1}{2}
\operatorname{Area}(\Omega)$. Hence by the values of the number
$\omega_K$ of roots of unity,
$$
\Vol((\pr{\Ga_K})_\H\bs\H)=\frac{\sqrt{|D_K|}}{2\,\omega_K}\;.
$$
Note that $\Vol(\pr G_\H\bs\H)=[(\pr{\Ga_K})_\H:\pr G_\H]
\Vol((\pr{\Ga_K})_\H\bs\H)$.

Let $\C=\C(f)$, which is indeed a real hyperbolic plane in $\htr$,
whose set of points at infinity is $\C_\infty(f)$ (hence $\infty$ is a
point at infinity of $\C(f)$ if and only if $a=0$). Note that $\C$
is invariant under the group $\autom$ by Equation \eqref{eq:antiequiv}
(which implies that $\C(f\circ g)= g^{-1}\C(f)$ for every
$g\in\SLOK$). The arithmetic group $\autom$ acts with finite covolume
on $\C(f)$, its finite subgroup $\{\pm\id\}$ acting trivially. By
definition, 
$$
\covol(\autom\cap G)= \Vol\big(\;_{\mbox{$\Pautom\cap F$}}
\bs\C(f)\big)\;.
$$  
Note that $\covol(\autom\cap G)$ depends only on the $G$-orbit of $f$,
by Equation \eqref{eq:antiequiv} and since $\operatorname{SU}_{f \circ
  g} (\OOO_K)=g^{-1}\autom g$ for every $g\in \SLOK$.  By its
definition, $R_G(f)$ is the index of the subgroup $\Pautom\cap \pr G$
in $\pr G_\C$, hence
$$
\Vol(F_\C\bs\C)=\frac{1}{R_G(f)}\covol(\autom\cap G)\;.
$$

\medskip Now that we have defined $n,k,F,\H $ and $\C$, let us pause
in the proof, by recalling the following easy exercice in group
theory.

\blemm\label{lemcountdoubleclass} Let $C$ be a group and let
$A,B,A',B'$ be subgroups of $C$, such that $A\subset A'$ and $B\subset
B'$, both with finite indices. Let $D$ be the set of elements $g\in C$
such that $g^{-1}A'g\cap B' =\{1\}$. Then the fibers of the canonical
map from $A\backslash D/ B$ to $A'\backslash D/ B'$ all have cardinal
$[A':A][B':B]$.  
\elemm

\dem Note that the subset $D$ of $C$, being invariant under left
translation by $A'$ and under right translation by $B'$, is a disjoint
union of double cosets $D=\coprod_{i\in I} A'g_iB'$. Write
$A'=\coprod_{j=1}^m Aa_j$ and $B'=\coprod_{k=1}^n b_kB$. Let us prove
that $D=\coprod_{i\in I,1\leq j\leq m,1\leq k\leq n} Aa_jg_ib_kB$,
which yields the result. It is clear that $D$ is the union of the
double cosets $Aa_jg_ib_kB$. Let us prove that for every $a\in A$ and
$b\in B$, if the equality $aa_jg_ib_kb=a_{j'}g_{i'}b_{k'}$ holds, then
$i=i',j=j',k=k'$, which implies the disjointness of these double
cosets. That equality implies first that $i=i'$ by the
definition of the double coset representatives $g_i$'s, and thus that
$g_i^{-1}a_{j'}^{-1}aa_jg_i=b_{k'}b^{-1}b_k^{-1}$. Since
$a_{j'}^{-1}aa_j$ and $b_{k'}b^{-1}b_k^{-1}$ belong to $A'$ and $B'$
respectively, the assumptions defining $D$ imply that they are both
equal to the identity element. Hence $aa_j=a_{j'}$ and
$b_kb=b_{k'}$. By the definition of the left representatives $a_j$'s
and the right representatives $b_k$'s, we hence have $j=j'$ and
$k=k'$. \cqfd

\medskip The last step of the proof of Theorem \ref{theo:mainversionG}
consists in relating the two counting functions $\psi_{f,G}$ and
$\N_{F,\C,\H}$, in order to apply Theorem \ref{theo:mainBPP}.  For
every $g\in\SLC$, let us compute the hyperbolic length of the common
perpendicular arc $\delta_{g^{-1},\id}$ between the real hyperbolic
plane $g^{-1}\C$ and the horosphere $\partial \H$, assuming that they
do not meet. The radius of the circle $\C_\infty(f\circ g)$, which is
the boundary at infinity of $g^{-1}\C$ by Equation
\eqref{eq:antiequiv}, is $\frac{\sqrt{\Delta(f\circ g)}} {|a(f\circ
  g)|}$, where $a(f\circ g)$ is the coefficient of $|u|^2$ in $f\circ
g(u,v)$. An immediate computation gives
\begin{equation}\label{eq:calcullongarcperpcomm}
\ell(\delta_{g^{-1},\id})=
\ln\frac{|a(f\circ g)|}{\sqrt{\Delta(f\circ g)}}=
\ln\frac{|f\circ g(1,0)|}{\sqrt{\Delta(f)}}\;.
\end{equation}
With the conventions that we have taken, this formula is also valid if
$g^{-1}\C$ and $\partial \H$ meet.

Let $U_K$ denote the stabilizer of $(1,0)\in\CC^2$ in $\GaK$, which
consists of the unipotent upper triangular matrices in $\GaK$. For
every $s> 0$, using Equation \eqref{eq:calcullongarcperpcomm}, we have
$$
\begin{aligned}
  \psi_{f,G}(s) &=\card\;\big\{[g]\in(\autom\cap G)\bs G/(U_K\cap
  G)\;:\;
  |f\circ g(1,0)|\leq s\big\}\\
  &= \card\;\big\{[g]\in(\Pautom\cap\pr G)\bs \pr G/(\pr{U_K}\cap\pr
  G)\;:\;
  \ell(\delta_{g^{-1},\id})\leq \ln\frac{s}{\sqrt{\Delta}}\big\}\;.
\end{aligned}
$$

Apply Lemma \ref{lemcountdoubleclass} to $C=\overline{G}=F$,
$A=\Pautom\cap\pr G$, $A'=F_\C$, $B=\pr{U_K}\cap\pr G$ and $B'=F_\H$.
Note that there are only finitely many elements $[g]\in F_\C\bs
F/F_\H$ such that the multiplicity $m(g^{-1},\id)=
1/{\card(\,g^{-1}F_\C g\; \cap\; F_\H})$ is different from $1$.
Therefore, by the above and by Theorem \ref{theo:mainBPP}, we have
$$
\!\!\!\!\!\begin{aligned}
\psi_{f,G}(s)& \sim R_G(f)\;[\pr G_\H:\pr{U_K}\cap \pr G]\;
  \card\;\big\{[g]\in F_\C\bs F/F_\H\;:\;
  \ell(\delta_{g^{-1},\id})\leq \ln\frac{s}{\sqrt{\Delta}}\big\}
\\ &
\sim R_G(f)\;[\pr G_\H:\pr{U_K}\cap \pr G]\;
  \N_{F,\C,\H}\big(\ln\frac{s}{\sqrt{\Delta}}\big)
\end{aligned}
$$
$$ 
\sim  R_G(f)\;[\pr G_\H:\pr{U_K}\cap \pr G]\;
\frac{2\;[(\pr{\Ga_K})_\H:\pr G_\H]\frac{\sqrt{|D_K|}}{2\;\omega_K}\;
\frac{1}{R_G(f)}\covol\big(\autom\cap G\big)}
{4\pi\;\frac{1}{\iota_G}[\GaK:G]\frac{1}{4\pi^2}|D_K|^{3/2} 
\;\zeta_K(2)}\;\;\frac{s^2}{\Delta},
$$
as $s$ tends to $+\infty$. Since $[(\pr{\Ga_K})_\H:\pr G_\H][\pr
G_\H:\pr{U_K}\cap \pr G]=\frac{1}{2}[\Ga_{K,\infty}:G_{(1,0)}]$,
Theorem \ref{theo:mainversionG} follows.  
\cqfd

\bigskip We now state the precise asymptotic result of the number of
nonequivalent representations of rational integers, satisfying some
congruence relations and having absolute value at most $s$, by a given
integral indefinite binary Hermitian form. Given a nonzero ideal
$\aaa$ in $\OOO_K$, let $\iota_\aaa=1$ if $2\in\aaa$, and $\iota_\aaa=2$
otherwise; consider the congruence subgroups
$$ 
\GaKa=\Big\{\begin{pmatrix}
  \alpha&\beta\\\gamma&\delta\end{pmatrix}\in\Ga_K:
\alpha-1,\delta-1,\gamma,\beta\in\mathfrak
a\Big\},\;\;\;\GaKaZ=\Big\{\begin{pmatrix}
  \alpha&\beta\\\gamma&\delta\end{pmatrix}\in\Ga_K: \gamma\in\mathfrak
a\Big\}\;.
$$ 
Both $\GaKaZ$ and $\GaKa$ coincide with $\GaK$ when $\aaa=\OOO_K$.

\bcoro \label{coro:maingeneral} Let $f$ be an integral indefinite
binary Hermitian form over an imaginary quadratic number field $K$,
and let $\aaa$ be a nonzero ideal in $\OOO_K$. As $s$ tends to
$+\infty$, we have
$$
\card\;\;_{\mbox{$\autom\cap\GaKa $}}\bs
\big\{(u,v)\in\P_K\;:\;u-1,v\in\aaa,\;\;|f(u,v)|\leq s\big\}$$
$$\sim 
\frac{\pi\;\iota_\aaa\;\covol(\autom\cap\GaKa)}
{2\;{\rm N}(\aaa)^2
\prod_{\ppp|\aaa}\big(1+\frac 1{{\rm N}(\ppp)}\big)\;
|D_K|\;\zeta_K(2)\;\Delta(f)}\;\;\;s^2\;,
$$
and 
$$
\card\;\;_{\mbox{$\autom\cap\GaKaZ$}}\bs
\big\{(u,v)\in\P_K\;:\;v\in \aaa,\;\;|f(u,v)|\leq s\big\}$$
$$\sim \frac{\pi\;\covol(\autom\cap\GaKaZ)}
{2\,{\rm N}(\aaa)
\prod_{\ppp|\aaa}\big(1-\frac 1{{\rm N}(\ppp)^2}\big)\;
|D_K|\;\zeta_K(2)\;\Delta(f)}\;\;\;s^2\;.
$$
\ecoro

\dem The orbits of $(1,0)\in\CC^2$ under $\GaKa$ and $\GaKaZ$ are
precisely the sets $\{(u,v)\in\P_K:u-1,v\in\aaa\}$ and
$\{(u,v)\in\P_K:v\in\aaa\}$ respectively.

The indices of $\GaKa$ and $\GaKaZ$ in $\GaK$, as computed for example
in Theorems VII.16 and VII.17 of \cite{Newman72}, are 
$$
[\GaK:\GaKa]={\rm N}(\aaa)^3
\prod_{\ppp|\aaa}\big(1+\frac 1{{\rm N}(\ppp)}\big)\;\;\;{\rm
and}\;\;\;
[\GaK:\GaKaZ]={\rm N}(\aaa)
\prod_{\ppp|\aaa}\big(1-\frac 1{{\rm N}(\ppp)^2}\big)\;,
$$
where the products are taken over the prime ideals
$\ppp$ of $\OOO_K$ dividing $\aaa$.

The index in $\Ga_{K,\infty}$ of the stabilizer of $(1,0)\in\CC^2$ in
$\GaKa$ is $\omega_K\,{\rm N}(\aaa)$. The index in $\Ga_{K,\infty}$ of
the stabilizer of $(1,0)\in\CC^2$ in $\GaKaZ$ is $\omega_K$. 

Note that $-\id$ belongs to $\GaKaZ$, but belongs to $\GaKa$ if and
only if $2\in\aaa$, so that $\iota_{\GaKaZ}=1$ and
$\iota_{\GaKa}=\iota_\aaa$.

The corollary now follows from Theorem \ref{theo:mainversionG},
applied with $G=\GaKa$ and $G=\GaKaZ$.  
\cqfd

\medskip Theorem \ref{theo:mainHerm} follows from Corollary
\ref{coro:maingeneral}, by taking in both results $\aaa=\OOO_K$.  Note
that from the techniques of \cite[Sect.~9]{ElsGruMen98}, only the much
weaker result $\psi_f(s)=O(s^2\log s)$ seems to be obtainable (see
\cite[Coro.~2.12]{ElsGruMen87}).

\bigskip In the following concluding remarks, for any positive integer
$\Delta$, let
$$
f_\Delta(u,v)=|u|^2-\Delta|v|^2,
$$
which is an integral indefinite binary Hermitian form with
discriminant $\Delta$.

\medskip
\noindent{\bf Remark 1.}  Here is a computation of $\covol(\autom)$
for $f$ an integral indefinite binary Hermitian form over $K$, with
discriminant $\Delta$, following \cite{Maclachlan84,MacRei03} instead
of \cite{Prasad89}.

Maclachlan has proved in \cite{Maclachlan84} that $\autom$ and
$\operatorname{SU}_{f_\Delta}(\OOO_K)$ are commensurable up to
conjugation, in the following way. Since the limit set of $\PSLOK$ is
$\phtr=\CC\cup\{\infty\}$ and since $\autom=\operatorname{SU}_{-f}
(\OOO_K)$, we may assume, up to replacing $f$ by an element in its
$\GaK$-orbit or its opposite, that $a=a(f)>0$.  Let $G_a$ be the
congruence subgroup of $\autom$ which is the preimage of the upper
triangular subgroup by the morphism $\autom\ra\operatorname{SL}_2
(\OOO_K/a\OOO_K)$ of reduction modulo $a$ of the coefficients. Let
$g=\Big(\begin{array}{cc} \frac{1}{\sqrt{a}} &
  -\frac{\overline{b}}{\sqrt{a}}\\0&\sqrt{a}\end{array}\Big)\in\SLC$.
By an easy computation, we have that $f\circ g= f_\Delta$, and that
$g^{-1} G_ag$ is contained in $\SLOK$. Hence $g^{-1} G_ag$ is a finite
index subgroup of $\operatorname{SU}_{f_\Delta}(\OOO_K)$. Therefore, we
have
$$
\covol(\autom)= \frac{[\operatorname{SU}_{f_\Delta}(\OOO_K):g^{-1}
  G_ag]}{[\autom:G_a]}\; \covol(\operatorname{SU}_{f_\Delta}(\OOO_K))\;.
$$
Maclachlan has also proved in \cite{Maclachlan84} that
$\operatorname{SU}_{f_\Delta}(\OOO_K)$ is commensurable with a lattice
derived from a quaternion algebra, in the following way. Let
$d_K=\frac{D_k}{4}$ if $D_K\equiv 0\mod 4$, and $d_K=D_K$
otherwise. Let $A$ be the quaternion algebra with Hilbert symbol
$\Big(\begin{array}{c}d_K,\Delta\\\hline \QQ\end{array}\Big)$ over
$\QQ$, which splits over $K$. Let $\Delta(A)$ be the reduced
discriminant of $A$. Let $\OOO$ be the order $\ZZ[i,j,ij]$ in $A$ for
the standard basis $1,i,j,ij$ of $A$, let $\OOO_{\rm max}$ be the
maximal order containing $\OOO$, and let $\OOO^1,\OOO^1_{\rm max}$ be
their groups of units. Let $\varphi:A\ra A\otimes_\QQ K=\M_2(K)$ be
the natural embedding, given by
$$
\alpha+\beta i+\gamma j+\delta ij\mapsto
\begin{pmatrix}\alpha+\beta\sqrt{d_K}&
  \Delta(\ga+\delta\sqrt{d_K}\,)\\\ga-\delta\sqrt{d_K}&
  \alpha-\beta\sqrt{d_K}\end{pmatrix}\,.
$$ 
An easy computation shows that $\varphi(\OOO^1)$ is a subgroup of
$\operatorname{SU}_{f_\Delta}(\OOO_K)$. Then by
\cite[Theo.~11.1.1]{MacRei03}, we have
$$
\covol(\operatorname{SU}_{f_\Delta}(\OOO_K))= 
\frac{\pi\;[\OOO^1_{\rm max}:\OOO^1]}
{3\;[\operatorname{PSU}_{f_\Delta}(\OOO_K):\pr{\varphi(\OOO^1)}\,]}\;
\prod_{p\,|\,\Delta(A)}(p-1)\;,
$$
where $p$ ranges over the positive rational primes. 
This gives a formula for $\covol(\autom)$.

\medskip

\noindent{\bf Remark 2.} 
Assume that $K=\QQ(i)$ in this remark.  Let us give in this case,
following \cite{MacRei91}, a more precise formula for $\covol(\autom)$
where $f$ is an integral indefinite binary Hermitian form over $K$,
with discriminant $\discr$. Combined with Theorem \ref{theo:mainHerm},
Corollary \ref{coro:gaussian} will follow. 

Since $\autom=\operatorname{SU}_{kf}(\OOO_K)$ for every $k\in\NN-\{0\}$,
we may assume, as required in \cite{MacRei91}, that $f$ is {\em
  primitive}, that is, with $f$ as in Equation \eqref{eq:form}, the
coefficients $a$, $c$, and the real and imaginary parts of $b$ have no
common divisor in the rational integers. Note that the subgroup
$\Pautom$ of $\pr{\GaK}$ is denoted 
by $\operatorname{Stab}(\C(f),\pr\GaK)$ in \cite[p.~161]{MacRei91}, and
it is a maximal Fuchsian subgroup of $\pr{\Ga_K}$ (loc.~cit.).

The hyperbolic plane $\C(f_\discr)$ associated to the form $f_\discr$
is the halfsphere of Euclidean radius $\sqrt{\discr}$ centered at
$0$. A formula due to Humbert (see for instance
\cite[p.~169]{MacRei91}) gives
\begin{equation}\label{eq:humbert}
\covol(\operatorname{SU}_{f_\discr}(\OOO_K))
=\eta\,\pi\,\discr
\prod_{\substack{p|\discr\\p\textrm{ odd}}}
\big(1+\bigg(\frac{-1}p\bigg)p^{-1}\big)\;,
\end{equation}
where $p$ ranges over the positive rational primes,
$\big(\frac{-1}p\big)$ is a Legendre symbol, $\eta=1/2$ if
$\discr\equiv 0\mod 4$ and $\eta=1$ otherwise.

The maximal Fuchsian subgroups of $\pr{\Ga_K}$ are classified in
\cite{MacRei91}, yielding the following cases.

If $\discr\equiv 0,3\mod 4$, then $\Pautom$ is a conjugate in
$\pr{\Ga_K}$ of $\operatorname{PSU}_{f_\discr}(\OOO_K)$, and its
covolume is given by Equation \eqref{eq:humbert} (see
\cite[p.170]{MacRei91}).

If $\discr\equiv 1\mod 4$, there are two cases: If the coefficients
$a$ and $c$ are even, then
$$
\covol(\autom) =
\frac 1 3 \covol(\operatorname{SU}_{f_\discr}(\OOO_K))\;;
$$
otherwise, $\Pautom$ is a conjugate in $\pr{\Ga_K}$ of
$\operatorname{PSU}_{f_\discr}(\OOO_K)$ (see \cite[p.~172]{MacRei91}). 

If $\discr\equiv 2\mod 4$, there are two cases: If the coefficients
$a$ and $c$ are even, then
$$
\covol(\autom)=
\frac 1{\eta'} \covol\operatorname{SU}_{f_\discr}(\OOO_K)\;,
$$
where $\eta'\in\{2,6\}$ satisfies $\eta'=\discr\mod 8$; otherwise,
$\Pautom$ is a conjugate in $\pr{\Ga_K}$ of
$\operatorname{PSU}_{f_\discr}(\OOO_K)$ (see \cite[p.~173]{MacRei91}).

This proves Corollary \ref{coro:gaussian} of the introduction.

{\small
\bibliography{../biblio}
}

\bigskip
{\small\noindent \begin{tabular}{l} 
Department of Mathematics and Statistics, P.O. Box 35\\ 
40014 University of Jyv\"askyl\"a, FINLAND.\\
{\it e-mail: parkkone@maths.jyu.fi}
\end{tabular}
\medskip

\noindent \begin{tabular}{l}
DMA, UMR 8553 CNRS\\
Ecole Normale Sup\'erieure, 45 rue d'Ulm\\
75230 PARIS Cedex 05, FRANCE\\
{\it e-mail: Frederic.Paulin@ens.fr}
\end{tabular} and \begin{tabular}{l}
D\'epartement de math\'ematique, B\^at.~425\\
Universit\'e Paris-Sud 11\\
91405 ORSAY Cedex, FRANCE\\
{\it e-mail: frederic.paulin@math.u-psud.fr}
\end{tabular}
}

\end{document}

%% file: macros_angl.tex
\usepackage{amssymb,epsfig}
\usepackage{amsmath}
\usepackage{comment}
\usepackage[dvips]{color}
\usepackage[T1]{fontenc}

\newtheorem{defi}{Definition}
\newtheorem{prop}[defi]{Proposition}
\newtheorem{theo}[defi]{Theorem}
\newtheorem{conj}[defi]{Conjecture}
\newtheorem{lemm}[defi]{Lemma}
\newtheorem{coro}[defi]{Corollary}
\newtheorem{rema}[defi]{Remark}
\newtheorem{exem}[defi]{Example}
\newtheorem{exems}[defi]{Examples}

\newcommand{\bdefi}{\begin{defi}}
\newcommand{\edefi}{\end{defi}}
\newcommand{\bprop}{\begin{prop}}
\newcommand{\eprop}{\end{prop}}
\newcommand{\btheo}{\begin{theo}}
\newcommand{\etheo}{\end{theo}}
\newcommand{\blemm}{\begin{lemm}}
\newcommand{\brema}{\begin{rema}}
\newcommand{\erema}{\end{rema}}
\newcommand{\bexer}{\begin{exem}}
\newcommand{\eexer}{\end{exem}}
\newcommand{\bexems}{\begin{exems}}
\newcommand{\eexems}{\end{exems}}
\newcommand{\bconj}{\begin{conj}}
\newcommand{\econj}{\end{conj}}
\newcommand{\elemm}{\end{lemm}}
\newcommand{\bcoro}{\begin{coro}}
\newcommand{\ecoro}{\end{coro}}
\newcommand{\dem}{\noindent{\bf Proof. }}

\usepackage{mathrsfs}
\renewcommand\mathcal{\mathscr}

\newcommand{\M}{{\cal M}}
\newcommand{\N}{{\cal N}}

\newcommand{\F}{{\cal F}}

\renewcommand{\H}{{\cal H}}
\newcommand{\OOO}{{\cal O}}
\newcommand{\C}{{\cal C}}

\renewcommand{\P}{{\cal P}}

\newcommand{\maths}[1]{{\mathbb #1}}  

\newcommand{\RR}{\maths{R}}
\newcommand{\NN}{\maths{N}}
\newcommand{\CC}{\maths{C}}
\newcommand{\QQ}{\maths{Q}}

\newcommand{\SSS}{\maths{S}}

\newcommand{\HH}{\maths{H}}

\newcommand{\ZZ}{\maths{Z}}
\newcommand{\PP}{\maths{P}}

\newcommand{\aaa}{{\mathfrak a}}

\newcommand{\ppp}{{\mathfrak p}}

\newcommand{\ra}{\rightarrow}
\newcommand{\bs}{\backslash}

\newcommand{\ov}[1]{{\overline #1}}

\newcommand{\ga}{\gamma}
\newcommand{\Ga}{\Gamma}

\newcommand{\card}{{\operatorname{Card}}}
\renewcommand{\Re}{{\operatorname{Re}}}

\newcommand{\Vol}{\operatorname{Vol}}

\newcommand{\covol}{\operatorname{Covol}}
\newcommand{\id}{\operatorname{id}}

\newcommand{\cqfd}{\hfill$\Box$}

\newcommand{\PSLC}{\operatorname{PSL}_{2}(\CC)}
\newcommand{\SLC}{\operatorname{SL}_{2}(\CC)}

\newcommand{\htr}{{\HH\,}^3_\RR}
\newcommand{\hnr}{{\HH\,}^n_\RR}

\newcommand{\phtr}{\partial\htr}

\newcommand{\SLOK}{\operatorname{SL}_{2}(\OOO_K)}
\newcommand{\PSLOK}{\operatorname{PSL}_{2}(\OOO_K)}

\newcommand{\autom}{\operatorname{SU}_f(\OOO_K)}
\newcommand{\Pautom}{\operatorname{PSU}_f(\OOO_K)}
\newcommand{\const}{\iota}
\newcommand{\discr}{\Delta}

\newcommand{\GaK}{{\Ga_{K}}}
\newcommand{\GaKa}{{\operatorname{\Ga_{\mathit K}(\mathfrak a)}}}
\newcommand{\GaKaZ}{{\operatorname{\Ga_{\mathit K,0}(\mathfrak a)}}}

\newcommand{\pr}{\overline}

\newcounter{fig}

